\documentclass[12pt,english]{amsart}
\usepackage[T1]{fontenc}
\usepackage[latin9]{inputenc}
\usepackage[a4paper]{geometry}
\usepackage{verbatim}
\usepackage{calc}
\usepackage{amsthm}
\usepackage{amssymb}
\usepackage{graphicx}
\usepackage{setspace}
\makeatletter

\usepackage[english]{babel}
\usepackage{amsfonts}
\usepackage{epigraph}

\usepackage{etoolbox}

\makeatletter
\patchcmd{\epigraph}{\@epitext{#1}}{\itshape\@epitext{#1}}{}{}
\makeatother

\title{Daniel Augusto da Silva, Poet of Mathematics}

\author[Carlos Florentino]{Carlos Florentino}

\address{Departamento de Matem\'{a}tica, Faculdade de Ci\^encias, Univ. de Lisboa, Campo Grande, Edf. C6, Lisbon, Portugal}

\email{caflorentino@ciencias.ulisboa.pt}

\thanks{This work was partially supported by CMAFcIO of the University of Lisbon, FCT, Portugal.}


\usepackage{babel}

\makeatother

\usepackage{babel}
\begin{document}
\begin{abstract}
Daniel da Silva was a remarkable Scientist and Mathematician of the
mid 19th century. Working in Portugal, isolated from the main scientific
centers of the time, his investigations in pure mathematics had almost
no impact. Apart from giving a short biography of his life and work,
this article makes the case for considering him an unavoidable character
in the History of Science and one of the founders of Discrete Mathematics,
through his introduction of a key method in Enumerative Combinatorics:
the Principle of Inclusion-Exclusion.
\end{abstract}

\dedicatory{Dedicated to the memory of Jaime de Lima Mascarenhas}
\maketitle
\begin{quote}
\epigraph{"Ninguem está authorisado a capitular quaesquer theorias mathematicas como destituidas de applicação vantajosa, como um mero recreio de elevadas intelligencias, e como inuteis trabalhos em relação à verdadeira sciencia."}{--- \textup{Daniel da Silva}}
\end{quote}

\section{Introduction}

Daniel Augusto da Silva (1814-1878) is one of the greatest Portuguese
scientists of the 19th Century. His contributions are remarkable,
for their quality and originality, not only in the context of Portuguese
Science, but also internationally. Even though his life and mathematical
work is documented by historians and some researchers, his important
contributions have not yet received the deserved recognition from
the mathematical community%
.

Daniel da Silva produced only a few manuscripts of scientific nature,
most of them published by the Lisbon Academy of Sciences between the
years 1851 and 1876. These articles belong to the fields of Statics
(more precisely, what we call today Geometric or Rational Mechanics),
Physics/Chemistry, Statistics and Actuarial Sciences, and Number Theory.

The book by Francisco Gomes Teixeira (1851\textendash 1933) \cite{T1},
published in 1934, and considered to be the most important reference
about History of Mathematics in Portugal up to the end of the 19th
century, singles out the four most important mathematical characters,
according to the author: Pedro Nunes (1502-1578), Anastácio da Cunha
(1744-1787), Monteiro da Rocha (1734-1819), and Daniel da Silva. Gomes
Teixeira refers to Daniel in this way: 

\emph{\textquotedblleft Daniel da Silva, poet of mathematics, searched
in these sciences what they have of beautiful; {[}\dots {]} he gave
to the world of numbers his Statics, without worrying with the applications
of this chapter of rational mechanics, that others later did, and
gave it also his beautiful investigations about binomial congruences.\textquotedblright }
\begin{figure}
\includegraphics[scale=0.48]{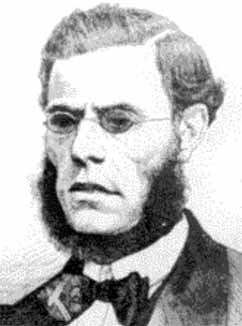}\caption{Daniel Augusto da Silva}
\end{figure}

As happened with other great names of Science, the work of Daniel
da Silva did not take place without unfortunate moments and drama.
Indeed, two facts played a major role: his enormous creative capabilities
ended up being constrained by a major illness he suffered for several
years; and the fact that he always wrote in Portuguese made it very
difficult for the recognition of his contributions in other European
countries with much more solid and developed scientific communities.

In this short article, we are going to summarise his biography and
his work, concentrating on his most important innovations in Pure
Mathematics, area which represented, in his own words, his great passion
(sections 2 and 3). For their relevance in today's mathematics, we
finish by analyzing with further detail (section 4) da Silva's contributions
to Discrete Mathematics and Number Theory, which deserve to be widely
known.

\subsection*{Acknowledgements}

I would like to thank many friends and colleagues who encouraged me
to write this article. I thank especially José Francisco Rodrigues,
Luís Saraiva, Helena Mascarenhas, Cristina Casimiro, Pedro J. Freitas
and Susana Ferreira for the careful review.

\section{Education and Professional Life}

Daniel Augusto da Silva was born in Lisbon on the 16th of May of 1814,
being the second son of Roberto José da Silva and Maria do Patrocínio.
Roberto Silva was a merchant although his specific business does not
seem to be documented.

Daniel's formative years took place against a difficult background
of tensions and wars in Portugal. In fact, historians agree that the
whole first half of the 19th century was not auspicious for the development
of scientific culture in the country: this period included the Napoleon
invasions, between 1807 and 1811, the violent campaigns opposing liberals
and absolutists, and a civil war between 1828 and 1834.

At the age of 15, Daniel enrolled in the Royal Navy Academy (Academia
Real da Marinha, ARM), and took mathematics courses ranging from Arithmetics
to Calculus, as well as some courses on Mechanics, Physics and Navigation.
He also attended courses at the Lisbon Royal Naval Observatory. He
immediately showed special talent for mathematics and was awarded
a distinction in each of the three years there. In 1832, he entered,
by merit, the Royal Academy of Marine Guards (Academia Real dos Guardas-Marinhas,
ARGM), an academy typically reserved for sons of officials, and was
appointed Navy Officer in 1833.
\begin{figure}
\includegraphics[scale=0.44]{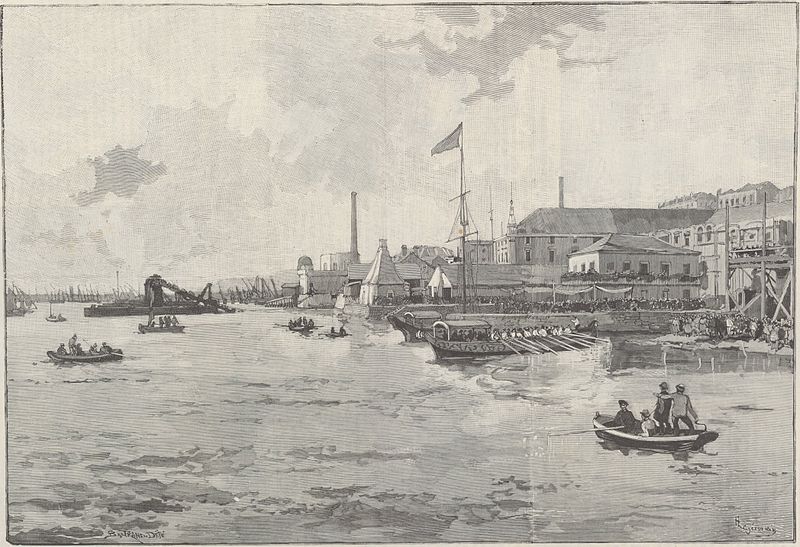}\caption{Dom Pedro II, Brasil's emperor, docks at Lisbon's Navy Arsenal in
1876.}
\end{figure}

As he became interested in mathematics, after finishing the ARGM degree
in 1835, he asked for permission and for a small fellowship, to enrol
in the Mathematics Faculty of the Coimbra University (the only Portuguese
University at the time). Being approved by the Navy, he moved to Coimbra,
and his performance in the University was no less brilliant than in
both Academies: many years later some of his old professors could
still remember the brightness of da Silva as a student.

Having finished his studies in Coimbra in 1839, he immediately returned
to Lisbon, and followed a career in the Navy. He was promoted to Brigadier
on 1840, and, later that year, to Second Lieutenant of the Navy. Following
the French trend of \textquotedblleft Grands Écoles\textquotedblright ,
in 1845, the ARGM was transformed into the Navy School (Escola Naval),
and Daniel was appointed as a professor there. He taught Mechanics;
Astronomy and Optics; Artillery and Fortification, and Geography and
Hydrography. Initially, he was hired as \textquotedblleft Lente Substituto\textquotedblright{}
(literally substitute teacher), and became \textquotedblleft Lente
Proprietario\textquotedblright{} (Full Professor) in 1848.

There were two unfortunate moments when da Silva missed the opportunity
of becoming a Professor at the blossoming Polytechnic School of Lisbon.
This School had just been created in 1837, in the context of a higher
education reform, by a Royal Decree, to replace the ARM. In 1839 he
applied to a teaching position through a competitive process but unfortunately,
for health reasons, Daniel da Silva could not be present in a kind
of interview/examination. Even after justifying his absence, the panel
decided to cancel the placement, afraid of a possible impugnation.
The second occasion was in 1848 when the Directing Body of the Polytechnic
School of Lisbon, acknowledging Daniel da Silva\textquoteright s value,
directly asked the Government to authorise his appointment, something
rejected on the grounds that, by law, all places should be filled
by public competition.

Nonetheless, while in the Navy School, it is between 1848 and 1852
that Daniel da Silva experiences his main creative period with the
completion of his 3 first manuscripts, sent for publication by the
Royal Lisbon Academy of Sciences (Academia Real das Ciências, ARC\footnote{Nowadays called ``Academia das Ciências de Lisboa''.}).
In 1851, he becomes a \emph{corresponding member} of this Society,
and is elected as \emph{full member} the following year. 

In late 1852 his health problems became so severe, and magnified by
his overwork and his great dedication to research, that he applied
for a leave and went to Madeira, hoping to recover there. However,
his poor health persisted and he was unable to carry out his duties
until the Naval Health Board classified him ``unfit for active duty''
in 1859. 

This same year, he was ellected \emph{honorary member} of the ARC,
and married Zefferina d'Aguiar (1825-1913) from the town of Funchal.
Daniel and Zefferina had a single child, Júlio Daniel da Silva who
was born in 1866. Sadly, Júlio would die at the age of 25 without
descendants.

Even without teaching duties, he continued to hold his Navy position
until retiring in 1868. In his latest years, worried that his passion
by Pure Mathematics would worsen his health condition, Daniel continued
to do research, but dedicated himself to more applied Sciences, publishing
works in actuarial sciences and on the theory of the flame. In his
words:

``\emph{The passion for the study of mathematics, that was in me
greatly disordered by excess, many years now has been reduced to the
modest proportions of a platonic love.}''

In 1871, the young Francisco Gomes Teixeira, a third year student
in Coimbra heard his professor José Queirós mention da Silva's theory
of couples in Statics with high praises, recalling his brightness
as a student, more than 30 years before. An excellent mathematician
himself, Teixeira became acquainted with Daniel's work, and decided
to write an essay on continued fractions, a subject of Daniel's incomplete
chapter 10 of \cite{dS}. He then wrote a letter to Daniel including
this essay, and this started an excellent and joyful relationship
between the two. Da Silva soon invited Gomes Teixeira to become a
member of the ARC and tried to get him a position in the Lisbon Astronomic
Observatory. After Daniel died, Teixeira presented his eulogy to the
ARC and he would become da Silva's most complete biographer \cite{T2}.
Further accounts on Daniel da Silva life and work can be found in
\cite{Di,Du,Ma,O,Va,Sa1}.

\section{Scientific Work}

In the period 1849-51, Daniel da Silva wrote 3 manuscripts concerning
investigations on geometric methods in statics of rigid bodies, and
on number theory. These show that he was an avid reader of the classics,
being inspired by names such as Euler (1707-83), Lagrange (1736-1813),
Legendre (1752-1833), Gauss (1777-1855) and Poinsot (1777-1859). He
would obtain mathematical articles published in European Academies
of Sciences, especially from the one in Paris.

Daniel's first paper \emph{``On the transformation and reduction
of binaries of forces''} was written before 1850, but only published
in 1856 by the ARC. Closely following an article of Louis Poinsot
on the same theme, this article contains no original results, but
presents a new approach including some simplified proofs. 

His second paper, \emph{``Memoir on the rotation of forces about
their points of application'', }was read to the Academy in 1850 and
published the following year. Here, da Silva considers a system of
forces turning around their points of application, but maintaining
their relative angles during the rotation. The article was written
without knowledge of the results of A. F. Möbius (1790\textendash 1868)
on this subject. Möbius had incorrectly stated that if a system is
in equilibrium in four different orientations, then it is in equilibrium
in all possible positions. In this memoir, da Silva describes correctly
the equilibrium properties of a system of forces, and also proves
that, in general, there are only four equilibrium positions.
\begin{figure}
\includegraphics[scale=0.3]{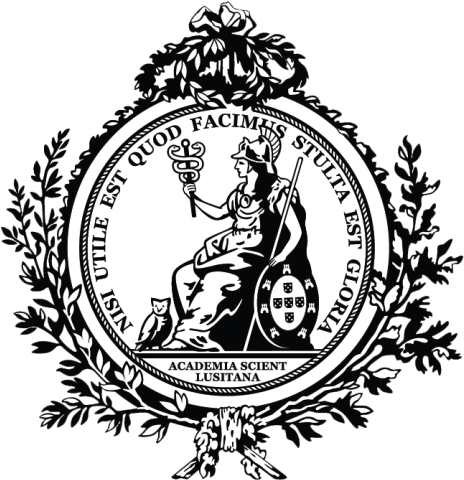}\caption{Logo of ``Academia das Ciências de Lisboa''}

\end{figure}

The third of da Silva's memoirs, on number theory, was read to the ARC
 on March 1852, but his illness prevented the completion of the published
version \cite{dS} (see below). We dedicate section 4 to an exposition
(of the initial part) of this article, and its important achievements.
About it, Teixeira writes: 

\emph{``The main subject he considered was the resolution of binomial
congruences, a theory which belongs simultaneously to the domain of
higher arithmetic and higher algebra, and he enriched it with such
important and general results that his name deserves to be included
in the list of those who founded it. It was indeed Daniel da Silva
who first gave a method to solve systems of linear congruences, an
honour which has been unduly attributed to the distinguished English
arithmetician Henry Smith}\footnote{\emph{Henry J. S. Smith (1826-1883), see \cite{Sm}.}}\emph{,
who only in 1861 dealt with this subject, and {[}Daniel da Silva{]}
was also the one who first undertook the general study of binomial
congruences.''\cite{T1}}

As mentioned, da Silva's health prevented him from doing any serious
mathematical research for quite some time. Only many years later,
his health did improve a little and he again undertook research. This
new period started in 1866, with a very short article on Statics and,
shortly after, two articles on Statistics and Actuarial Sciences.
Daniel proposed mathematical models of demography and applied them
to the financial structuring of pension funds, and in particular to
one of the oldest portuguese Welfare Institutions called Montepio
Geral (founded in 1840). These two articles, \emph{``Average annual
amortization of pensions in the main Portuguese Welfare Institute}s'',
and\emph{ ``Contribution to the comparative study of population dynamics
in Portugal'' }are described in detail in a recent PhD thesis \cite{Ma},
which is also a very important source of information on Daniel da
Silva, and on Actuarial Calculus and the Navy Schools, at the epoch. 

In 1872, he published \emph{``On several new formulae of Analytic
Geometry relative to oblique coordinate axes'', }which generalizes
certain well known formulae in Analytic Geometry to a setup based
on non-orthogonal frames.

He also carried out studies in the area of Physics and Chemistry,
that could have been motivated by his previous lectures in the Navy
School. He performed several experiments, with the collaboration of
António Augusto Aguiar (1838-1887), professor at the Polithecnic School
of Lisbon, and studied the speed of transmission of a gas flame in
its blueish and brightest part in the 1873 article\emph{ ``Considerations
and experiments about the flame''}.

In 1877, the last year of his life, Daniel received some unpleasant
news from France. As he had worked in scientific isolation, and used
the portuguese language, J. G. Darboux (1842-1917) had just published
results very similar to Daniel's own research on Statics (including
the same correction of Möbius' mistake) without acknowledging da Silva's
work. His letter to Teixeira denotes his disappointment: 

\emph{``Almost all propositions of Darboux are published twenty six
years ago in the Memoirs of the Lisbon Academy of Sciences, in my
work on the rotation of forces about their points of application!
{[}\dots {]} My memoir, which contains many other things, besides
those considered by Möbius, including a correction of one mistake
he did, the same one which Darboux proudly claims correcting, lies
ignored, for nearly twenty six years, in the Libraries of almost all
Academies of the world! What worth it is writing in portuguese!''}

Daniel da Silva knew of Darboux\textquoteright s article \cite{Da}
in the French journal ``Les Mondes'', and sent to it a ``Réclamation
de Priorité''; this reclamation was published in this periodical
on March that same year, but, as L. Saraiva writes in \cite{Sa1}:
\emph{``It certainly would have been better if he had written directly
to the French Academy of Sciences, where his work could have been
more widely discussed.``}

After Daniel's death, this whole story had such a profound influence
on Gomes Teixeira, that he became one of the first portuguese scientists
and mathematicians to continuously promote the interaction of portuguese
academics with foreign researchers. He founded the first mathematical
journal, independent of any academic institution, printed in the Iberian
Peninsula (Jornal de Sciencias Mathematicas e Astronomicas, see \cite{R})
which substantially contributed to the dissemination of portuguese
research. And he also stimulated the analysis of Daniel's work in
the international scene, inviting some of his students and collaborators
to review and to continue da Silva's work.

For a comparative analysis of the results of da Silva on Statics and
those of A. F. Möbius and F. Minding (1806\textendash 1885), see F.
A. Vasconcelos \cite{Va}, who was encouraged by Gomes Teixeira to
perform this detailed study.

Daniel's research on the propagation of the flame was rediscovered
by the German chemist and professor in Zurich, Karl Heumann (1850-1894)
who, recongnizing the priority of Daniel da Silva in some of these
investigations, wrote him a letter in 1878. Unfortunately, Daniel
would never read it, as the letter arrived right after he had passed
away. For a complete list of da Silva's publications, see \cite{DdS}.

\section{Discrete Mathematics and Number Theory}

The memoir entitled \emph{``General properties and direct resolution
of binomial congruences''} \cite{dS} was presented to the Lisbon
Academy of Sciences in 1852 and published two years later. There are
many reasons to consider this as Daniel da Silva's masterpiece.
\begin{figure}
\framebox[0.62\columnwidth]{\includegraphics[scale=0.65]{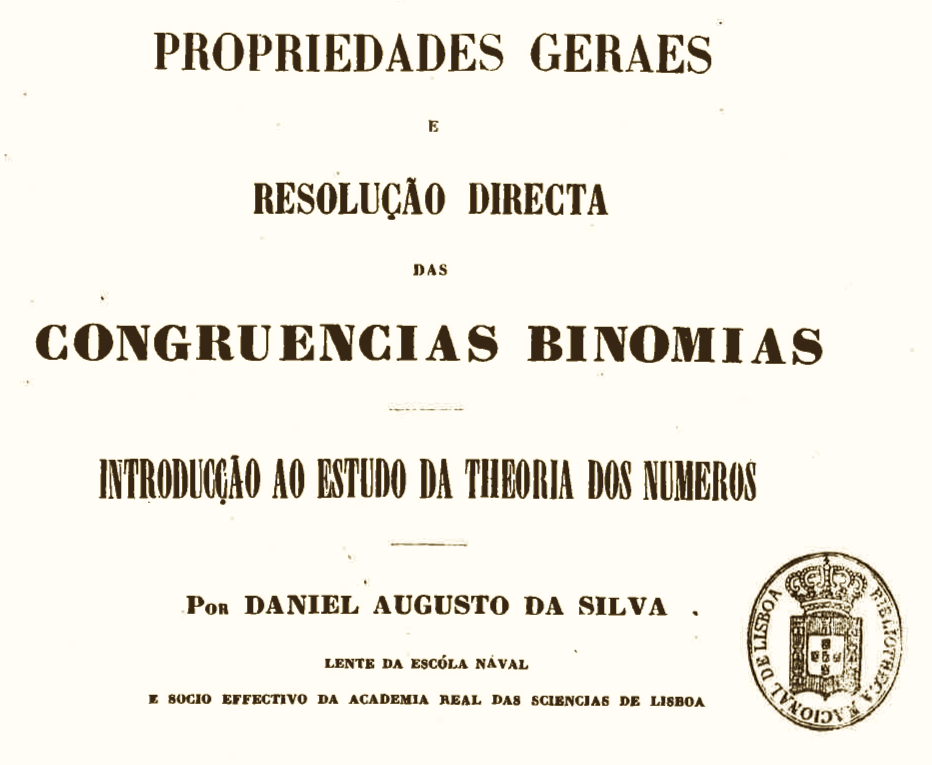}}\caption{Cover page of Daniel da Silva's memoir on Congruences}
\end{figure}

Right in the first pages, da Silva develops a case for the importance
of Pure Mathematics and its relationship with Applications, for the
relevance and elegance of Number Theory, citing and praising several
famous mathematicians such as Fermat (1601-65), Euler, Lagrange, Legendre,
Poinsot and Gauss, instead of going directly to the results as in
his other articles. Even though it includes many original results,
this memoir appears also to have a pedagogical goal, as indicated
by the subtitle \emph{``Introduction to the study of number theory\textquotedblright .
}This (as well as the expression \emph{``General Properties''})
hints that Daniel da Silva had in mind the foundations of a whole
new theory of mathematics, naturally abstract, but that could provide,
in his opinion, numerous applications in many contexts.

A second reason is that a big portion of the basics of what we call
today Discrete Mathematics are literally present in this work, constituting
a great advance for the epoch. A typical syllabus for a freshman Discrete
Mathematics course includes:
\begin{itemize}
\item Some Logic and Set Theory (including operations with sets (intersections,
unions, etc), cardinality, examples such as $\mathbb{Z}$, $\mathbb{Q}$,
$\mathbb{Z}/n\mathbb{Z}$ etc);
\item Basics of Number Theory (including modular arithmetic, divisibility,
the Euler $\varphi$ function, the theorems of Fermat/Euler, etc); 
\item Some Enumerative Combinatorics (including the binomial formula, the
principle of inclusion-exclusion, generating functions, etc).
\end{itemize}
It is remarkable that da Silva\textquoteright s memoir from 1854 treats,
in a clear, elegant and modern way, most of the above list of topics
and subtopics. Moreover, some of the methods used bear a striking
coincidence with those of textbooks for first year courses of Discrete/Finite
Mathematics adopted nowadays in Colleges and Universities around the
world. This happens, e.g, in Daniel's proof of the Euler\textquoteright s
formula for $\varphi(n)$, the function that counts the number of
positive integers less than $n\in\mathbb{N}$ and prime with it (see
below).

A final argument for considering \cite{dS} as Daniel's masterpiece
are the original results and their current relevance. Indeed, two
very important results introduced here are unanimously attributed
to da Silva: the famous Principle of Inclusion-Exclusion, and a formula
for congruences that generalises the well-known formula of Euler involving
his $\varphi$ function. Let us recall this material, from a modern
perspective, and compare with the way Daniel introduces it.

\subsection{The Principle of Inclusion-Exclusion}

The Principle of Inclusion-Exclusion (abbreviated PIE) is one of the
essential counting methods in Combinatorics, allowing a multitude
of applications. It is ubiquitous in most books dedicated to this
area.\footnote{One standard textbook \cite{St} explicitly mentions the PIE more
that 50 times. } The PIE generalises the formula for the cardinality of the union
of two finite sets $A$ and $B$:
\[
|A\cup B|=|A|+|B|-|A\cap B|,
\]
where $|A|$ denotes the cardinality of $A$. This is also widely
used in Probability or Measure Theory since, with appropriate interpretation,
we can replace cardinality by probability or by measure. To state
the PIE in modern terms, for $i=1,\cdots,N$ ($N\in\mathbb{N}$) consider
finite subsets $A_{1},\cdots,A_{N}$ of a given finite set $X$, and
let $A=\cup_{i=1}^{N}A_{i}$ be their union inside $X$, then the
PIE is the formula:
\begin{equation}
|A|=\sum_{i}|A_{j}|-\sum_{i<j}|A_{ij}|+\sum_{i<j<k}|A_{ijk}|-\cdots+(-1)^{N-1}|A_{12\cdots N}|\label{eq:PIE}
\end{equation}
where $A_{ij}:=A_{i}\cap A_{j}$, $A_{ijk}:=A_{i}\cap A_{j}\cap A_{k}$,
etc. Equivalently, the PIE determines the cardinality of the complement
$A^{c}:=X\setminus A$ of $A$, as: 
\begin{equation}
|A^{c}|=|X|-|A|=|X|-\sum_{i}|A_{i}|+\sum_{i<j}|A_{ij}|-\cdots+(-1)^{N}|A_{12\cdots N}|\label{eq:PIEc}
\end{equation}
For us, the most interesting aspect of da Silva's proof of \eqref{eq:PIE}
is that it requires the introduction of the \emph{notion of set},
at least of a finite one, a couple of decades before the foundations
of set theory laid by George Cantor (1845-1918)! Even though Daniel
believes the concept of set is of major importance, he refers to this
just as a convenient ``notation''. In da Silva's own words:

\emph{``To prove this formula we will employ a notation, that may
advantageously serve in other cases. Suppose that in a sequence $S$
of numbers (that we consider }united\emph{ and not }summed\emph{,
since if even some of them could be negative, there wouldn't result
any }subtraction\emph{) one asks which are the ones that satisfy some
property $a$; we will denote by $S_{a}$ the }reunion\emph{ of those
numbers; Similarly $S_{b}$, $S_{b,c}$, $S_{a,b,c}$ etc, the reunion
of those terms of $S$ verifying property $b$, or simultaneously
the properties $b$, $c$, etc}''.\footnote{We reproduced here Daniel's own emphasis in the 4 words: \emph{united},
\emph{summed}, \emph{subtraction} and \emph{reunion}.}

Then, using also the notation $^{\cdots,c,b,a}S$ for the elements
of $S$ that do not possess any of the properties $a,b,c,...$, da
Silva's presents his \textbf{symbolic formula} as follows:
\begin{equation}
^{\cdots,c,b,a}S=S[1-_{a}][1-_{b}][1-_{c}]\cdots\label{eq:Daniel-symbolic}
\end{equation}
and clarifies that, on the rigth hand side, \emph{``in such a product
the letters $a,\,b$ etc, become indices, and any composed index such
as $a_{b_{c}}$ becomes a simple index $a,b,c$'' }since\emph{ ``it
is easy to see that $S_{a_{b_{c}}}=S_{a,b,c}$''} and similarly for
upper left indices. 

Daniel continues, using $\psi$ for the cardinality of a set:

``\emph{The same formula also immediately gives us the number of
numbers contained in $^{\cdots,c,b,a}S$; denoting this number by
$\psi{}^{\cdots,c,b,a}S$, and letting the sign $\psi$ have an analogous
meaning applied to the additive and subtractive series in \eqref{eq:Daniel-symbolic}
it is clear that we will have:
\begin{equation}
\psi{}^{\cdots,c,b,a}S=\psi S[1-_{a}][1-_{b}][1-_{c}]\cdots".\label{eq:PIE-Daniel}
\end{equation}
}One of the most common modern proofs of PIE uses the characteristic
function of a subset $A\subset X$, defined by:
\[
\chi_{A}(x):=\begin{cases}
1, & x\in A\\
0, & x\notin A.
\end{cases}
\]
It is easy to see that characteristic functions are multiplicative
for intersections, and additive for disjoint unions:
\[
\chi_{A\cap B}=\chi_{A}\chi_{B},\quad\chi_{A\sqcup B}=\chi_{A}+\chi_{B}.
\]
This implies, upon observing that $A^{c}=\cap_{i=1}^{N}A_{i}^{c}$,
when $A=\cup_{i=1}^{N}A_{i}$:
\begin{equation}
\chi_{A^{c}}=\prod_{i=1}^{N}\chi_{A_{i}^{c}}=\prod_{i=1}^{N}(\mathbf{1}-\chi_{A_{i}})=\mathbf{1}-\sum_{i}\chi_{A_{i}}+\sum_{i<j}\chi_{A_{ij}}-\cdots+(-1)^{N}\chi_{A_{12\cdots N}},\label{eq:char-func}
\end{equation}
where $\mathbf{1}=\chi_{X}$ (the constant function 1 defined on $X$).
Noting that cardinality is just given by summing over $X$ (recall
this is a finite set): 
\[
|A|=\sum_{x\in X}\chi_{A}(x),
\]
equation \eqref{eq:char-func} is tranformed into the PIE, in the
form \eqref{eq:PIEc}. By replacing $a$ with $\chi_{A_{1}}$, $b$
with $\chi_{A_{2}}$ etc, the similarity of this proof with Daniel's
formulae \eqref{eq:Daniel-symbolic} and \eqref{eq:PIE-Daniel} is
manifest!
\begin{figure}
\framebox[0.71\columnwidth]{\includegraphics[scale=1.16]{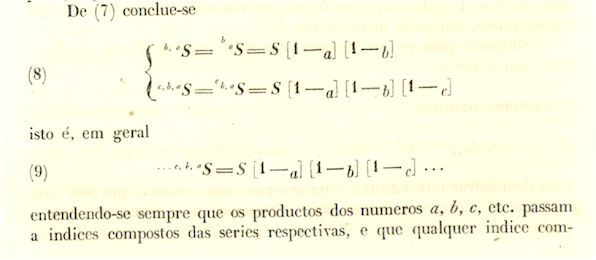}}\caption{Part of page 10 of \cite{dS} where Daniel's symbolic formula is written }
\end{figure}

Next, Daniel applies his formula \eqref{eq:PIE-Daniel} to deduce
the formula of Euler
\begin{equation}
\varphi(n)=n\left(1-\frac{1}{a}\right)\left(1-\frac{1}{b}\right)\cdots\left(1-\frac{1}{c}\right),\label{eq:phi-Euler}
\end{equation}
when $n=a^{\alpha}b^{\beta}\cdots c^{\gamma}$ is the prime factorization
of $n\in\mathbb{N}$, and this is done in \emph{the same way} as in
many modern books: see, for example, the standard textbook \cite{Bi},
from 2005, that derives the formula \eqref{eq:phi-Euler} precisely
as da Silva does, using PIE.

\subsection{Euler's Theorem and Bézout identity}

Euler used the function $\varphi$ in his celebrated theorem:
\begin{equation}
a^{\varphi(n)}\equiv1\quad(\mbox{mod }n),\label{eq:Euler}
\end{equation}
for $a\in\mathbb{Z}$ relatively prime to $n\in\mathbb{N}$. This
formula generalizes Fermat's theorem:
\[
a^{p-1}\equiv1\quad(\mbox{mod }p),
\]
for $p$ a prime and $a$ not multiple of $p$. 

The following elegant and interesting generalization of \eqref{eq:Euler}
was proved by Daniel da Silva. Let $n=a_{1}a_{2}\cdots a_{k}$ where
all factors ($k>1$ in number) are pairwise relatively prime. Then:
\begin{equation}
a_{1}^{\varphi(n/a_{1})}+a_{2}^{\varphi(n/a_{2})}+\cdots+a_{k}^{\varphi(n/a_{k})}\equiv k-1\quad(\mbox{mod }n).\label{eq:Daniel-Euler}
\end{equation}
The proof can be found in the original memoir of Daniel da Silva \cite{dS}
or in \cite[p. 211]{CRS}.

The particular case $k=2$ of this formula connects beautifully with
the Bézout identity. This identity (in a simplified form) states that,
given two relatively prime natural numbers $a,b\in\mathbb{N}$, there
are solutions $x,y\in\mathbb{Z}$ to: 
\[
ax+by=1,
\]
and it is well known that this can be solved by an ancestral method:
the (extended) Euclidean algorithm. From \eqref{eq:Daniel-Euler},
we see that Daniel's formula provides a direct solution: 
\[
x=a^{\varphi(b)-1},\quad y=b^{\varphi(a)-1},
\]
for the ``congruence version'' of Bézout's identity:
\[
ax+by\equiv1\quad(\mbox{mod }ab)
\]
(under the same assumptions on $a,b$). 

The monograph goes on with many interesting applications of these
formulas and related questions. Among these, Daniel provides direct
resolutions for linear congruences, for the chinese remainder theorem,
and for many congruences of the form 
\[
ax^{n}\equiv b\quad(\mbox{mod }N).
\]
Daniel's health problems intensified as he was approaching the end
of his monograph: he wasn't able to revise neither the preface nor
the final part (see \cite{Sa2}). For the same reason, the last two
chapters are incomplete: in the 9th some theorems he would like to
add are missing; and sections 4 and 5 of the 10th chapter have only
their title. The last section was supposed to be a study of continued
fractions, the theme of Teixeira's first letter to Daniel, 20 years
later.

Without any impact whatsoever at the time it was written, this memoir
was discovered almost by accident, half a century later, by the italian
mathematician Cristoforo Alasia (1864-1918) who, suprised by its depth,
dedicated three articles to Daniel's work between 1903 and 1914 (all
in italian, the first being \cite{A}). However, as far as we know,
only one conference proceedings (written in portuguese) has addressed
the concept of set in da Silva's work \cite{dC}.\footnote{I thank J. F. Rodrigues for the indication of this reference.}

We finish with a quote from \cite{dS}, a wonderful illustration of
Daniel da Silva's passion for mathematics and his opinion on the importance
of pure mathematics and its role in science:

``\emph{In general, it can be said that nobody is authorized to capitulate
any mathematical theory as deprived from advantageous applications,
as a mere recreation of elevated minds, and as useless work towards
true science. All acquired truths are that many elements of acumulated
intelectual wealth. Soon or late, the day will come when concrete
science will have to search within this vast arsenal the necessary
tools for great discoveries, which in this way will pass from speculative
theorems to the category of practical truths. Every day, one or another
area of mathematical-physics or celestial or industrial mechanics,
is observed to suddenly stop its development to implore assistance
from the further improvements in pure analysis, without which those
very important sciences cannot progress.}''

Daniel da Silva could have hardly guessed that this phrase would be,
more than a century later, so appropriate to the subject of his own
article. In fact, the algorithms that preserve the security of data
across the internet, such as the famous RSA (Rivest\textendash Shamir\textendash Adleman)
criptographic system that we use (even without noticing) on a daily
basis, depend crucially on Euler's formula \eqref{eq:Euler}.%

\end{document}